\input amstex
\documentstyle{amsppt}
%----------------------------------------------------------------
% Title:     On two elliptic curves associated with perfect cuboids.
% Author:    Ruslan Sharipov
% Comments:  AmSTeX, 11 pages, amsppt style
% MSC-class: 11D25, 11D72, 12E05, 14E05, 14G05, 14H52
%----------------------------------------------------------------
%           Replacement for output macro definition
%
\catcode`@=11
\redefine\output@{%
  \def\break{\penalty-\@M}\let\par\endgraf
  \ifodd\pageno\global\hoffset=105pt\else\global\hoffset=8pt\fi  
  \shipout\vbox{%
    \ifplain@
      \let\makeheadline\relax \let\makefootline\relax
    \else
      \iffirstpage@ \global\firstpage@false
        \let\rightheadline\frheadline
        \let\leftheadline\flheadline
      \else
        \ifrunheads@ %\let\makefootline\relax
        \else \let\makeheadline\relax
        \fi
      \fi
    \fi
    \makeheadline \pagebody \makefootline}%
  \advancepageno \ifnum\outputpenalty>-\@MM\else\dosupereject\fi
}
\def\Beta{\mathchar"0\hexnumber@\rmfam 42}
\catcode`\@=\active
%----------------------------------------------------------------
\nopagenumbers
\chardef\textvolna='176
\def\negskp{\hskip -2pt}

\chardef\degree="5E
%\def\compos{\,\raise 1pt\hbox{$\sssize\circ$} \,}
%\def\id{\operatorname{id}}
%\font\eightrm=cmr8
%\def\LT{\operatorname{\text{\eightrm LT}}}
%\def\LM{\operatorname{\text{\eightrm LM}}}
%\def\LC{\operatorname{\text{\eightrm LC}}}
%\accentedsymbol\hatgamma{\kern 2pt\hat{\kern -2pt\gamma}}
%\accentedsymbol\checkgamma{\kern 2.5pt\check{\kern -2.5pt\gamma}}
\def\blue#1{#1}

\catcode`#=11\def\diez{#}\catcode`#=6
\catcode`&=11\catcode`&=4
\catcode`_=11\def\podcherkivanie{_}\catcode`_=8
\catcode`~=11\def\volna{~}\catcode`~=\active
\def\mycite#1{\cite{\blue{#1}}\immediate\special{ps:
     ShrHPSdict begin /ShrBORDERthickness 0 def}}
\def\myciterange#1#2#3#4{\cite{\blue{#2#3#4}}\immediate\special{ps:
     ShrHPSdict begin /ShrBORDERthickness 0 def}}
\def\mytag#1{%
    \tag#1}
\def\mythetag#1{\thetag{\blue{#1}}\immediate\special{ps:
     ShrHPSdict begin /ShrBORDERthickness 0 def}}
\def\myrefno#1{\no#1}
\def\myhref#1#2{\blue{#2}\immediate\special{ps:
     ShrHPSdict begin /ShrBORDERthickness 0 def}}
\def\myEarXivlink{\myhref{http://arXiv.org}{http:/\negskp/arXiv.org}}

\def\mytheorem#1{\csname proclaim\endcsname{Theorem #1}}
\def\mytheoremwithtitle#1#2{\csname proclaim\endcsname{Theorem #1#2}}
\def\mythetheorem#1{\blue{#1}\immediate\special{ps:
     ShrHPSdict begin /ShrBORDERthickness 0 def}}
\def\mylemma#1{\csname proclaim\endcsname{Lemma #1}}
\def\mylemmawithtitle#1#2{\csname proclaim\endcsname{Lemma #1#2}}

\def\mycorollary#1{\csname proclaim\endcsname{Corollary #1}}

\def\mydefinition#1{\definition{Definition #1}}

\def\myconjecture#1{\csname proclaim\endcsname{Conjecture #1}}
\def\myconjecturewithtitle#1#2{\csname proclaim\endcsname{Conjecture #1#2}}

\def\myproblem#1{\csname proclaim\endcsname{Problem #1}}
\def\myproblemwithtitle#1#2{\csname proclaim\endcsname{Problem #1#2}}
\def\mytheproblem#1{\blue{#1}\immediate\special{ps:
     ShrHPSdict begin /ShrBORDERthickness 0 def}}

%----------------------------------------------------------------
% Cyrillic fonts definition
\font\eightcyr=wncyr8
%\font\tencyr=wncyr10
%----------------------------------------------------------------
\pagewidth{360pt}
\pageheight{606pt}
\topmatter
\title
On two elliptic curves associated with perfect cuboids.
\endtitle
\rightheadtext{On two elliptic curves associated with perfect cuboids.}
\author
Ruslan Sharipov
\endauthor
\address Bashkir State University, 32 Zaki Validi street, 450074 Ufa, Russia
\endaddress
\email\myhref{mailto:r-sharipov\@mail.ru}{r-sharipov\@mail.ru}
\endemail
\abstract
     A rational perfect cuboid is a rectangular parallelepiped whose edges and face 
diagonals are given by rational numbers and whose space diagonal is equal to unity. 
Finding such a cuboid is equivalent to finding a perfect cuboid with all integer 
edges and diagonals, which is an old unsolved problem. Recently, based on a symmetry 
approach, it was shown that edges and face diagonals of rational perfect cuboid are
roots of two cubic equations whose coefficients depend on two rational parameters. 
Six special cases where these cubic equations are reducible have been already found. 
Two more possible cases of reducibility for these cubic equations are considered
in the present paper. They lead to a pair of elliptic curves.
\endabstract
\subjclassyear{2000}
\subjclass 11D25, 11D72, 12E05, 14E05, 14G05, 14H52\endsubjclass
\endtopmatter
%\loadbold
%\loadeufb
\TagsOnRight
\document
% \input countstyle

%\special{header=resource.eps}
\head
1. Introduction.
\endhead
     The problem of a perfect cuboid is known since 1719. For the history of
cuboid studies the reader is referred to \myciterange{1}{1}{--}{44}. Let 
$x_1$, $x_2$, $x_3$ be edges of a cuboid and $d_1$, $d_2$, $d_3$ be its face 
diagonals. Recently, as a result of the series of papers \myciterange{45}{45}{--}{50},
two cubic equations for $x_1$, $x_2$, $x_3$ and $d_1$, $d_2$, $d_3$ were derived:
$$
\align
&\hskip -2em
x^3-E_{10}\,x^2+E_{20}\,x-E_{30}=0,
\mytag{1.1}\\
&\hskip -2em
d^{\kern 1pt 3}-E_{01}\,d^{\kern 1pt 2}+E_{02}\,d-E_{03}=0.
\mytag{1.2}
\endalign
$$
The numbers $x_1$, $x_2$, $x_3$ are roots of the equation \mythetag{1.1}, while
$d_1$, $d_2$, $d_3$ are roots of the equation \mythetag{1.2}. Apart from \mythetag{1.1}
and \mythetag{1.2}, the numbers $x_1$, $x_2$, $x_3$ and $d_1$, $d_2$, $d_3$ should obey 
the following auxiliary equations:
$$
\hskip -2em
\aligned
&x_1\,x_2\,d_3+x_2\,x_3\,d_1+x_3\,x_1\,d_2=E_{21},\\
&x_1\,d_2+d_1\,x_2+x_2\,d_3+d_2\,x_3+x_3\,d_1+d_3\,x_1=E_{11},\\
&x_1\,d_2\,d_3+x_2\,d_3\,d_1+x_3\,d_1\,d_2=E_{12}.
\endaligned
\mytag{1.3}
$$
The left hand sides of the equations \mythetag{1.3} are three of nine elementary 
multisymmetric polynomials that correspond to the permutation group $S_3$ acting 
upon the numbers $x_1$, $x_2$, $x_3$ and $d_1$, $d_2$, $d_3$ broken into pairs:
$(x_1,d_1)$, $(x_2,d_2)$, $(x_3,d_3)$. For the theory of general  multisymmetric 
polynomials the reader is referred to \myciterange{51}{51}{--}{71}.\par
      The coefficients of the cubic equations \mythetag{1.1} and \mythetag{1.2}
and the right hand sides of the equations \mythetag{1.3} are rational functions 
of two rational parameters $b$ and $c$. \pagebreak They are given by explicit 
formulas. Here are the formulas for $E_{01}$, $E_{10}$, and $E_{11}$:
$$
\align
&\hskip -2em
E_{11}=-\frac{b\,(c^2+2-4\,c)}{b^2\,c^2+2\,b^2-3\,b^2\,c+c-b\,c^2\,+2\,b},
\mytag{1.4}\\
\vspace{2ex}
&\hskip -2em
E_{01}=-\frac{b\,(c^2+2-2\,c)}{b^2\,c^2+2\,b^2-3\,b^2\,c+c-b\,c^2+2\,b},
\mytag{1.5}\\
\vspace{2ex}
&\hskip -2em
E_{10}=-\frac{b^2\,c^2+2\,b^2-3\,b^2\,c\,-c}{b^2\,c^2+2\,b^2-3\,b^2\,c
+c-b\,c^2+2\,b}.
\mytag{1.6}
\endalign
$$
The formulas \mythetag{1.4}, \mythetag{1.5}, and \mythetag{1.6}, were derived in
\mycite{49} by solving the equation\footnotemark\
$$
\hskip -2em
(2\,E_{11})^2+(E_{01}^2+L^2-E_{10}^2)^2-8\,E_{01}^2\,L^2=0, 
\mytag{1.7}
$$
\footnotetext{\ In deriving \mythetag{1.4}, \mythetag{1.5}, \mythetag{1.6} the
parameter $L$ in the equation \mythetag{1.7} was taken for the unity. However,
the corresponding solution for the general case $L\neq 1$ easily follows from
\mythetag{1.4}, \mythetag{1.5}, \mythetag{1.6} by homogeneity.}
which  was derived in \mycite{48}. Below are the formulas for $E_{02}$, 
$E_{21}$, $E_{12}$:
$$
\pagebreak
\gather
\gathered
E_{02}=\frac{1}{2}\,(28\,b^2\,c^2-16\,b^2\,c-2\,c^2-4\,b^2-b^2\,c^4
+\,4\,b^3\,c^4-12\,b^3\,c^3\,+\\
+\,4\,b\,c^3+24\,b^3\,c-8\,b\,c-2\,b^4\,c^4
+12\,b^4\,c^3-26\,b^4\,c^2-8\,b^2\,c^3\,+\\
+\,24\,b^4\,c-16\,b^3-8\,b^4)\,(b\,c-1-b)^{-2}\,(b\,c-c-2\,b)^{-2},
\endgathered\qquad\quad
\mytag{1.8}\\
\vspace{1ex}
\gathered
E_{21}=\frac{b}{2}\,(5\,c^6\,b-2\,c^6\,b^2+52\,c^5\,b^2-16\,c^5\,b
-2\,c^7\,b^2+2\,b^4\,c^8\,+\\
+\,142\,b^4\,c^6-26\,b^4\,c^7-426\,b^4\,c^5-61\,b^3\,c^6+100\,b^3\,c^5
+14\,c^7\,b^3\,-\\
-\,c^8\,b^3-20\,b\,c^2-8\,b^2\,c^2-16\,b^2\,c-128\,b^2\,c^4-200\,b^3\,c^3\,+\\
+\,244\,b^3\,c^2+32\,b\,c^3-112\,b^3\,c+768\,b^4\,c^4-852\,b^4\,c^3
+568\,b^4\,c^2\,+\\
+\,104\,b^2\,c^3-208\,b^4\,c+8\,c^4-4\,c^3+16\,b^3+32\,b^4-2\,c^5)\,\times\\
\times\,(b^2\,c^4-6\,b^2\,c^3+13\,b^2\,c^2-12\,b^2\,c+4\,b^2+c^2)^{-1}\,\times\\
\times\,(b\,c-1-b)^{-2}\,(b\,c-c-2\,b)^{-2},
\endgathered\qquad\quad
\mytag{1.9}\\
\vspace{2ex}
\hskip -2em
\gathered
E_{12}=(16\,b^6+32\,b^5-6\,c^5\,b^2+2\,c^5\,b-62\,b^5\,c^6+62\,b^6\,c^6\,-\\
-\,180\,b^6\,c^5+18\,b^5\,c^7-12\,b^6\,c^7-2\,b^5\,c^8+b^6\,c^8+248\,b^5\,c^2\,+\\
+\,248\,b^6\,c^2-96\,b^6\,c+321\,b^6\,c^4-180\,b^5\,c^3-144\,b^5\,c
-360\,b^6\,c^3\,+\\
+\,b^4\,c^8+8\,b^4\,c^6-6\,b^4\,c^7+18\,b^4\,c^5+7\,b^3\,c^6+90\,b^5\,c^5
-14\,b^3\,c^5\,-\\
-\,c^7\,b^3+17\,b^2\,c^4+28\,b^3\,c^3-28\,b^3\,c^2-4\,b\,c^3+8\,b^3\,c
-57\,b^4\,c^4\,+\\
+\,36\,b^4\,c^3+32\,b^4\,c^2-12\,b^2\,c^3-48\,b^4\,c-c^4+16\,b^4)\,\times\\
\times\,(b^2\,c^4-6\,b^2\,c^3+13\,b^2\,c^2-12\,b^2\,c+4\,b^2+c^2)^{-1}\,\times\\
\times\,(b\,c-1-b)^{-2}\,(b\,c-c-2\,b)^{-2}.
\endgathered\qquad
\mytag{1.10}\\
\endgather
$$
The formulas \mythetag{1.8}, \mythetag{1.9}, \mythetag{1.10} were derived in 
\mycite{50} by substituting \mythetag{1.4}, \mythetag{1.5}, \mythetag{1.6} 
and the length of the space diagonal $L=1$ into the appropriate formulas from 
\mycite{48}. The formulas for $E_{20}$, $E_{30}$, and $E_{03}$ are similar:
$$
\gather
\gathered
E_{20}=\frac{b}{2}\,(b\,c^2-2\,c-2\,b)\,(2\,b\,c^2-c^2-6\,b\,c+2
+4\,b)\,\times\\
\times\,(b\,c-1-b)^{-2}\,(b\,c-c-2\,b)^{-2},
\endgathered\qquad\quad
\mytag{1.11}\\
\vspace{1ex}
\gathered
E_{30}=c\,b^2\,(1-c)\,(c-2)\,(b\,c^2-4\,b\,c+2+4\,b)\,\times\\
\times\,(2\,b\,c^2-c^2-4\,b\,c+2\,b)\,\times\\
\times\,(b^2\,c^4-6\,b^2\,c^3+13\,b^2\,c^2-12\,b^2\,c+4\,b^2+c^2)^{-1}\,\times\\
\times\,(b\,c-1-b)^{-2}\,(-c+b\,c-2\,b)^{-2},
\endgathered\qquad\quad
\mytag{1.12}\\
\vspace{1ex}
\gathered
E_{03}=\frac{b}{2}\,(b^2\,c^4-5\,b^2\,c^3+10\,b^2\,c^2-10\,b^2\,c+4\,b^2+2\,b\,c\,+\\
+\,2\,c^2-b\,c^3)\,(2\,b^2\,c^4-12\,b^2\,c^3+26\,b^2\,c^2-24\,b^2\,c\,+\\
+\,8\,b^2-c^4\,b+3\,b\,c^3-6\,b\,c+4\,b+c^3-2\,c^2+2\,c)\,\times\\
\times\,((b^2\,c^4-6\,b^2\,c^3+13\,b^2\,c^2-12\,b^2\,c+4\,b^2+c^2)^{-1}\,\times\\
\times\,(b\,c-1-b)^{-2}\,(-c+b\,c-2\,b)^{-2}.
\endgathered\qquad\quad
\mytag{1.13}\\
\endgather
$$
The formulas \mythetag{1.11}, \mythetag{1.12}, \mythetag{1.13} were also derived in 
\mycite{50} by substituting \mythetag{1.4}, \mythetag{1.5}, \mythetag{1.6} into the 
appropriate formulas from \mycite{48}.\par
     Based on the formulas \mythetag{1.4} through \mythetag{1.6} and \mythetag{1.8}
through \mythetag{1.13}, the following two problems were formulated. 
\myproblem{1.1} Find all pairs of rational numbers $b$ and $c$ for which the
cubic equations \mythetag{1.1} and \mythetag{1.2} with the coefficients given
by the formulas \mythetag{1.6}, \mythetag{1.11},	\mythetag{1.12} and \mythetag{1.5}, 
\mythetag{1.8},	\mythetag{1.13} possess positive rational roots $x_1$, $x_2$, 
$x_3$, $d_1$, $d_2$, $d_3$ obeying the auxiliary polynomial equations 
\mythetag{1.3} whose right hand sides are given by the formulas \mythetag{1.9}, 
\mythetag{1.4}, \mythetag{1.10}. 
\endproclaim
\myproblem{1.2} Find at least one pair of rational numbers $b$ and $c$ for which 
the cubic equations \mythetag{1.1} and \mythetag{1.2} with the coefficients given
by the formulas \mythetag{1.6}, \mythetag{1.11},	\mythetag{1.12} and \mythetag{1.5}, 
\mythetag{1.8},	\mythetag{1.13} possess positive rational roots $x_1$, $x_2$, 
$x_3$, $d_1$, $d_2$, $d_3$ obeying the auxiliary polynomial equations 
\mythetag{1.3} whose right hand sides are given by the formulas \mythetag{1.9}, 
\mythetag{1.4}, \mythetag{1.10}. 
\endproclaim
     The problems~\mytheproblem{1.1} and \mytheproblem{1.2} are equivalent to the
appropriate problems for perfect cuboids. Therefore they are equally difficult. 
However, now we can consider more simple problems, e\.\,g\. we can search for the
cases where the polynomials in the left hand sides of the equations \mythetag{1.1}
and \mythetag{1.2} are reducible. 
\mydefinition{1.1} A polynomial with rational coefficients is called reducible over
$\Bbb Q$ if it splits into a product of two or more polynomials with rational 
coefficients.
\enddefinition
     The polynomials \mythetag{1.1} and \mythetag{1.2} can be reducible simultaneously
or in a separate way. Six simple cases of reducibility were found in \mycite{72}. In 
each of these six cases the polynomials \mythetag{1.1} and \mythetag{1.2} are reducible 
simultaneously. Unfortunately, or maybe fortunately, since otherwise the problem would be 
closed, \pagebreak none of them leads to a perfect cuboid. Here is the list of reducibility
conditions for all of these cases:
$$
\xalignat 2
&\hskip -2em
\text{1) \ }b=0\text{\ and }c\neq 0;
&&\text{4) \ }c=2\text{\ and }b\neq 1;\\
&\hskip -2em
\text{2) \ }c=0\text{\ and }b\,(1+b)\neq 0;
&&\text{5) \ }b\,(c-2)^2=-2\text{\ and }c\neq 2;
\mytag{1.14}\\
&\hskip -2em
\text{3) \ }c=1\text{\ and }b\neq -1;
&&\text{6) \ }2\,b\,(c-1)^2=c^2\text{\ and }c\,(c-1)\neq 0.
\qquad
\endxalignat
$$
In addition to \mythetag{1.14}, there are two more options that could lead to 
reducibility of the second polynomial \mythetag{1.2}. They are considered in 
this paper. As for the first polynomial \mythetag{1.1}, until the converse is 
proved, the reducibility of \mythetag{1.2} does not imply the reducibility 
of \mythetag{1.1} in general.\par
\head
2. Elliptic reducibility curves.
\endhead
     Let's consider the formula \mythetag{1.13} for the last term $E_{03}$ in the 
of the polynomial \mythetag{1.2}. If $E_{03}=0$, then the polynomial \mythetag{1.2}
has the rational root $d=0$, i\.\,e\. it is reducible over $\Bbb Q$. Looking
at \mythetag{1.13}, we see that the numerator in this formula is the product of three
terms. One of them is $b$. The vanishing condition $b=0$ for $E_{03}$ is already listed
in \mythetag{1.14}. Let's consider the other two vanishing conditions:
$$
\gather
\hskip -2em
b^2\,c^4-5\,b^2\,c^3+10\,b^2\,c^2-10\,b^2\,c+4\,b^2-b\,c^3+2\,b\,c+2\,c^2=0,
\mytag{2.1}\\
\vspace{2ex}
\hskip -2em
\gathered
2\,b^2\,c^4-12\,b^2\,c^3+26\,b^2\,c^2-24\,b^2\,c\,+\\
+\,8\,b^2-b\,c^4+3\,b\,c^3-6\,b\,c+4\,b+c^3-2\,c^2+2\,c=0.
\endgathered
\mytag{2.2}
\endgather
$$
The equations \mythetag{2.1} and \mythetag{2.2} are quadratic with respect to $b$.
Their discriminants with respect to $b$ are given by the following formulas:
$$
\align
&\hskip -2em
D_7=-(7\,c^4-40\,c^3+84\,c^2-80\,c+28)\,c^2,
\mytag{2.3}\\
&\hskip -2em
D_8=(c^4-8\,c^3+12\,c^2-16\,c+4)\,(c-1)^2\,(c-2)^2.
\mytag{2.4}
\endalign
$$
The indices $7$ and $8$ in \mythetag{2.3} and \mythetag{2.4} mean that we consider
the seventh and the eighth reducibility cases continuing the list \mythetag{1.14}.
\par
     Let's denote through $P_7$ and $P_8$ the square free factors of the discriminants
$D_7$ and $D_8$ respectively. Both of them are fourth order polynomials of $c$:
$$
\align
&\hskip -2em
P_7(c)=-7\,c^4+40\,c^3-84\,c^2+80\,c-28,
\mytag{2.5}\\
&\hskip -2em
P_8(c)=c^4-8\,c^3+12\,c^2-16\,c+4.
\mytag{2.6}
\endalign
$$
Using the polynomials \mythetag{2.5} and \mythetag{2.6}, we can write the polynomial
equations 
$$
\align
&\hskip -2em
y^2=-7\,c^4+40\,c^3-84\,c^2+80\,c-28,
\mytag{2.7}\\
&\hskip -2em
y^2=c^4-8\,c^3+12\,c^2-16\,c+4,
\mytag{2.8}
\endalign
$$
which are more simple than the equations \mythetag{2.1} and \mythetag{2.2}. The 
equations \mythetag{2.5} and \mythetag{2.6} define two elliptic curves on the 
$(y,c)$ plane (see \mycite{73}). The discriminants of the quartic polynomials 
$P_7(c)$ and $P_8(c)$ do coincide and are nonzero:
$$
\hskip -2em
D(P_7)=D(P_8)=-1048576=-2^{\kern 1pt 20}\neq 0.
\mytag{2.9}
$$
Due to \mythetag{2.9} the elliptic curves defined by the equations \mythetag{2.7}
and \mythetag{2.8} both are non-degenerate.\par
     In Number Theory an elliptic curve is understood as a genus $1$ curve with
at least one rational point either finite or at infinity (see \mycite{73}). In the 
case of the curve \mythetag{2.8} such a rational point is obvious:
$$
\xalignat 2
&\hskip -2em
y=2,
&&c=0.
\mytag{2.10}
\endxalignat
$$
In the case of the curve \mythetag{2.7} it is not obvious, but it does exist:
$$
\xalignat 2
&\hskip -2em
y=1,
&&c=1.
\mytag{2.11}
\endxalignat
$$
\mytheorem{2.1} Each rational point $(y,c)$ of the curve \mythetag{2.7} with
$c\neq 1$ and $c\neq 2$ produces a rational solution $(b,c)$ for the equation
\mythetag{2.1}, where 
$$
\hskip -2em
b=\frac{c\,(c^2+y-2)}{2\,(c-1)\,(c-2)\,((c-1)^2+1)}.
\mytag{2.12}
$$
\endproclaim
\mytheorem{2.2} Each rational point $(y,c)$ of the curve \mythetag{2.8} with
$c\neq 1$ and $c\neq 2$ produces a rational solution $(b,c)$ for the equation
\mythetag{2.2}, where 
$$
\hskip -2em
b=\frac{c^2+y-2}{4\,(c-2)\,(c-1)}.
\mytag{2.13}
$$
\endproclaim
     The proof for both theorems~\mythetheorem{2.1} and \mythetheorem{2.2} is
pure calculations. Note that the equalities \mythetag{2.12} and \mythetag{2.13} 
in them are linear with respect to $y$. Resolving these equalities for $y$, we 
can formulate the following converse theorems. 
\mytheorem{2.3} Each rational solution $(b,c)$ of the equation \mythetag{2.1} with
$c\neq 0$ produces a rational point $(y,c)$ for the elliptic curve \mythetag{2.7}, 
where 
$$
\hskip -2em
y=\frac{(2\,c^4-10\,c^3+20\,c^2-20\,c+8)\,b-c^3+2\,c}{c}.
\mytag{2.14}
$$
\endproclaim
\mytheorem{2.4} Each rational solution $(b,c)$ of the equation \mythetag{2.2} produces 
a rational point $(y,c)$ for the elliptic curve \mythetag{2.8}, where 
$$
\hskip -2em
y=(4\,c^2-12\,c+8)\,b-c^2+2.
\mytag{2.15}
$$
\endproclaim
The formulas \mythetag{2.12}, \mythetag{2.13}, \mythetag{2.14}, and \mythetag{2.15}
mean that the curve \mythetag{2.1} is birationally equivalent to the elliptic curve
\mythetag{2.7}, while the curve \mythetag{2.2} is birationally equivalent to the 
elliptic curve \mythetag{2.8}.\par
\head
3. Exceptional solutions and points.
\endhead
     Let's return to the cubic equations \mythetag{1.1} and \mythetag{1.2} and to
the auxiliary equations \mythetag{1.3}. We can call them cuboid equations since
they were derived from the original cuboid equations through the symmetry factorization
procedure as a result of the series of papers \myciterange{45}{45}{--}{50}. \pagebreak
The simultaneous non-vanishing condition for all denominators in the formulas 
\mythetag{1.4} through \mythetag{1.6} and in the formulas \mythetag{1.8} 
through \mythetag{1.13} is written as the following inequality:
$$
\hskip -2em
\gathered
(b^2\,c^4-6\,b^2\,c^3+13\,b^2\,c^2-12\,b^2\,c+4\,b^2+c^2)\,(b\,c-1-b)\,\times\\
\times\,(b\,c-c-2\,b)\,(b^2\,c^2+2\,b^2-3\,b^2\,c+c-b\,c^2+2\,b)\neq 0.
\endgathered
\mytag{3.1}
$$
Combining the inequality \mythetag{3.1} with the equation \mythetag{2.1}, we find the 
only rational exceptional point on the curve \mythetag{2.1}. It is the origin:
$$
\xalignat 2
&\hskip -2em
b=0,
&&c=0.
\mytag{3.2}
\endxalignat
$$
Similarly, combining \mythetag{3.1} with the equation \mythetag{2.2}, we find the only
rational exceptional point on the curve \mythetag{2.2}. It coincides with the point
\mythetag{3.2}.\par
      If $c=0$, the equation \mythetag{2.1} has the only rational solution with $b=0$,
i\.\,e\. this solution coincides with \mythetag{3.2}. It contradicts the inequality
\mythetag{3.1}.\par
      For $c=1$ and $c=2$ the equation \mythetag{2.1} has the following rational solutions:
$$
\xalignat 2
&\hskip -2em
b=-2,
&&c=1,
\mytag{3.3}\\
&\hskip -2em
b=2,
&&c=2.
\mytag{3.4}
\endxalignat
$$
The solutions \mythetag{3.3} and \mythetag{3.4} do not contradict the inequality
\mythetag{3.1}. But they are covered by the cases 3 and 4 listed in \mythetag{1.14}.
Therefore we call them exceptional solutions of the equation \mythetag{2.1} or
exceptional points of the curve \mythetag{2.1}.\par
      The formula \mythetag{2.14} maps the solutions \mythetag{3.3} and \mythetag{3.4}
of the equation \mythetag{2.1} to the following rational points of the elliptic curve 
\mythetag{2.7}:
$$
\xalignat 2
&\hskip -2em
y=1,
&&c=1,
\mytag{3.5}\\
&\hskip -2em
y=-2,
&&c=2.
\mytag{3.6}
\endxalignat
$$
The point \mythetag{3.5} coincides with \mythetag{2.11}. Along with \mythetag{3.5} 
and \mythetag{3.6}, the elliptic curve \mythetag{2.7} has the following two solutions 
being mirror images of the previous two:
$$
\xalignat 2
&\hskip -2em
y=-1,
&&c=1,
\mytag{3.7}\\
&\hskip -2em
y=2,
&&c=2.
\mytag{3.8}
\endxalignat
$$
Note that the rational points \mythetag{3.5}, \mythetag{3.6}, \mythetag{3.7}, and 
\mythetag{3.8} are exceptional in the sense of the theorem~\mythetheorem{2.1}. And
finally, note that the elliptic curve \mythetag{2.7} has no rational points
with $c=0$. Then we can formulate the following result.
\mytheorem{3.1} Non-exceptional rational solutions of the equation \mythetag{2.1},
if they exist, are in one-to-one correspondence with non-exceptional rational 
points of the curve \mythetag{2.7}. The correspondence is established by the 
formulas \mythetag{2.14} and \mythetag{2.12}. 
\endproclaim
     Let's proceed to the equation \mythetag{2.2} associated with the curve 
\mythetag{2.8}. Along with the solution \mythetag{3.2}, it has the following 
rational solution with $c=0$:
$$
\xalignat 2
&\hskip -2em
b=-\frac{1}{2},
&&c=0.
\mytag{3.9}
\endxalignat
$$
The solution \mythetag{3.9} does not contradict the inequality \mythetag{3.1}. 
But it is covered by the case 2 listed in \mythetag{1.14}. Therefore we call 
it an exceptional solution of the equation \mythetag{2.2} or an exceptional 
point of the curve \mythetag{2.2}.\par
      The formula \mythetag{2.15} maps the solution \mythetag{3.9} of the 
equation \mythetag{2.2} to the following rational point of the elliptic curve 
\mythetag{2.8}:
$$
\xalignat 2
&\hskip -2em
y=-2,
&&c=0,
\mytag{3.10}
\endxalignat
$$
Along with \mythetag{3.10}, the elliptic curve \mythetag{2.8} has the following 
solution being a mirror image of \mythetag{3.10} and coinciding with 
\mythetag{2.10}:
$$
\xalignat 2
&\hskip -2em
y=2,
&&c=0.
\mytag{3.11}
\endxalignat
$$
The points \mythetag{3.10} and \mythetag{3.11} of the curve \mythetag{2.8} are not
exceptional in the sense of the theorem~\mythetheorem{2.2}. The formula \mythetag{2.13}
maps them to the solutions \mythetag{3.9} and \mythetag{3.2} of the equation
\mythetag{2.2} respectively. But the latter ones are exceptional. Therefore we call 
the points \mythetag{3.10} and \mythetag{3.11} exceptional by convention. Note also 
that the the equations \mythetag{2.2} has no solutions with $c=1$ or $c=2$ and 
the curve \mythetag{2.8} has no rational points with $c=1$ or $c=2$ as well. Therefore
we can formulate the following theorem similar to the theorem~\mythetheorem{3.1}.
\mytheorem{3.2} Non-exceptional rational solutions of the equation \mythetag{2.2},
if they exist, are in one-to-one correspondence with non-exceptional rational 
points of the curve \mythetag{2.8}. The correspondence is established by the 
formulas \mythetag{2.15} and \mythetag{2.13}. 
\endproclaim
\head
4. The seventh reducibility case. 
\endhead
     The seventh and the eighth reducibility cases considered below are based on
non-exceptional rational solutions of the equations \mythetag{2.1} and \mythetag{2.2} 
or, which is equivalent, on non-exceptional rational points of the elliptic curves 
\mythetag{2.7} and \mythetag{2.8}. We do not study the problem of existence for such
solutions and/or such points in the present paper. Therefore the results below are 
conditional provided these point and these solutions do exist.\par
     Like in \mycite{72}, let's denote through $P(x)$ and $Q(d)$ the cubic polynomials 
in the left hand sides of the equations \mythetag{1.1} and \mythetag{1.2} respectively.
Then due to \mythetag{1.6}, \mythetag{1.11}, and \mythetag{1.12} the coefficients of
the polynomial $P(x)$ are functions of $b$ and $c$. Due to \mythetag{1.5}, 
\mythetag{1.8}, and \mythetag{1.13} the coefficients of $Q(d)$ are also functions of
$b$ and $c$.\par 
     The seventh reducibility case occurs if $(b,c)$ is a non-exceptional rational 
solution of the equation \mythetag{2.1}. This solution is produced from some 
non-exceptional rational point $(y,c)$ on the elliptic curve \mythetag{2.7} by means
of the formula \mythetag{2.12}. Let's substitute \mythetag{2.12} into the coefficients
of the polynomials $P(x)$ and $Q(d)$ and take into account the equation \mythetag{2.7}
in calculating $y^2$, $y^3$, $y^4$, $y^5$, etc. As a result we get two very huge 
expressions for $P(x)$ and $Q(d)$, but they turn out to be factorable in $x$ and $d$
so that we can formulate the following theorem.
\mytheorem{4.1} If the rational numbers $b$ and $c$ are produced from some 
non-exceptional rational point $(y,c)$ of the elliptic curve \mythetag{2.7} by means 
of the formula \mythetag{2.12}, then the cubic polynomials in \mythetag{1.1} and 
\mythetag{1.2} are reducible over $\Bbb Q$ and the equations \mythetag{1.1} and 
\mythetag{1.2} are factored as
$$
\xalignat 2
&\hskip -2em
P_{\kern 1pt 7.2}(x)\,(x+1)=0,
&&Q_{7.2}(d)\,\,d=0.
\mytag{4.1}
\endxalignat
$$
\endproclaim
     The second formula \mythetag{4.1} is not surprising since the equation
\mythetag{2.1}, and hence the equation \mythetag{2.7}, were derived from the condition
$E_{03}=0$. As for the first formula \mythetag{4.1}, it is proved by substituting
$x=-1$ into the formula for $P(x)$ transformed by means of the formula \mythetag{2.12} 
as described above. Selecting those rational points of the curve \mythetag{2.7} where
the polynomials $P(x)=P_{\kern 1pt 7.2}(x)\,(x+1)$ and $Q(d)=Q_{7.2}(d)\,\,d$ split into 
three linear factors is a separate problem. We do not consider this problem in the 
present paper.\par
\head
5. The eighth reducibility case. 
\endhead
     The eighth reducibility case occurs if $(b,c)$ is a non-exceptional rational 
solution of the equation \mythetag{2.2}. This solution is produced from some 
non-exceptional rational point $(y,c)$ on the elliptic curve \mythetag{2.8} by means
of the formula \mythetag{2.13}. Let's substitute \mythetag{2.13} into the coefficients
of the polynomials $P(x)$ and $Q(d)$ and take into account the equation \mythetag{2.8}
in calculating $y^2$, $y^3$, $y^4$, $y^5$, etc. As a result we get two very huge 
expressions for $P(x)$ and $Q(d)$, but they turn out to be factorable in $x$ and $d$ 
so that we can formulate the following theorem.
\mytheorem{5.1} If the rational numbers $b$ and $c$ are produced from some 
non-exceptional rational point $(y,c)$ of the elliptic curve \mythetag{2.8} by means 
of the formula \mythetag{2.13}, then the cubic polynomials in \mythetag{1.1} and 
\mythetag{1.2} are reducible over $\Bbb Q$ and the equations \mythetag{1.1} and 
\mythetag{1.2} are factored as
$$
\xalignat 2
&\hskip -2em
P_{8.2}(x)\,(x-1)=0,
&&Q_{8.2}(d)\,\,d=0.
\mytag{5.1}
\endxalignat
$$
\endproclaim
     The second formula \mythetag{5.1} is not surprising since the equation
\mythetag{2.2}, and hence the equation \mythetag{2.8}, were derived from the condition
$E_{03}=0$. As for the first formula \mythetag{5.1}, it is proved by substituting
$x=1$ into the formula for $P(x)$ transformed by means of the formula \mythetag{2.13} 
as described above. Selecting those non-exceptional rational points of the curve 
\mythetag{2.8} where the polynomials $P(x)=P_{\kern 1pt 8.2}(x)\,(x-1)$ and 
$Q(d)=Q_{8.2}(d)\,\,d$ split into three linear factors is a separate problem. We do not 
consider this problem in the present paper.\par
\head
9. Concluding remarks. 
\endhead
     The theory of rational points on elliptic curves is a very advanced and still 
developing area in modern mathematics. It comprises some intriguing open questions,
e\.\,g\. the Birch and Swinnerton-Dyer conjecture, which is one of the seven 
Millennium Prize Problems (see \mycite{74} and \mycite{75}). The author expects
that the observations and results of the present paper will make the perfect cuboid 
problem closer to this fascinating area of mathematics. 


\Refs
\ref\myrefno{1}\paper
\myhref{http://en.wikipedia.org/wiki/Euler\podcherkivanie 
brick}{Euler brick}\jour Wikipedia\publ 
Wikimedia Foundation Inc.\publaddr San Francisco, USA 
\endref
\ref\myrefno{2}\by Halcke~P.\book Deliciae mathematicae oder mathematisches 
Sinnen-Confect\publ N.~Sauer\publaddr Hamburg, Germany\yr 1719
\endref
\ref\myrefno{3}\by Saunderson~N.\book Elements of algebra, {\rm Vol. 2}\publ
Cambridge Univ\. Press\publaddr Cambridge\yr 1740 
\endref
\ref\myrefno{4}\by Euler~L.\book Vollst\"andige Anleitung zur Algebra, \rm
3 Theile\publ Kaiserliche Akademie der Wissenschaf\-ten\publaddr St\.~Petersburg
\yr 1770-1771
\endref
\ref\myrefno{5}\by Pocklington~H.~C.\paper Some Diophantine impossibilities
\jour Proc. Cambridge Phil\. Soc\. \vol 17\yr 1912\pages 108--121
\endref
\ref\myrefno{6}\by Dickson~L.~E\book History of the theory of numbers, 
{\rm Vol\. 2}: Diophantine analysis\publ Dover\publaddr New York\yr 2005
\endref
\ref\myrefno{7}\by Kraitchik~M.\paper On certain rational cuboids
\jour Scripta Math\.\vol 11\yr 1945\pages 317--326
\endref
\ref\myrefno{8}\by Kraitchik~M.\book Th\'eorie des Nombres,
{\rm Tome 3}, Analyse Diophantine et application aux cuboides 
rationelles \publ Gauthier-Villars\publaddr Paris\yr 1947
\endref
\ref\myrefno{9}\by Kraitchik~M.\paper Sur les cuboides rationelles
\jour Proc\. Int\. Congr\. Math\.\vol 2\yr 1954\publaddr Amsterdam
\pages 33--34
\endref
\ref\myrefno{10}\by Bromhead~T.~B.\paper On square sums of squares
\jour Math\. Gazette\vol 44\issue 349\yr 1960\pages 219--220
\endref
\ref\myrefno{11}\by Lal~M., Blundon~W.~J.\paper Solutions of the 
Diophantine equations $x^2+y^2=l^2$, $y^2+z^2=m^2$, $z^2+x^2
=n^2$\jour Math\. Comp\.\vol 20\yr 1966\pages 144--147
\endref
\ref\myrefno{12}\by Spohn~W.~G.\paper On the integral cuboid\jour Amer\. 
Math\. Monthly\vol 79\issue 1\pages 57-59\yr 1972 
\endref
\ref\myrefno{13}\by Spohn~W.~G.\paper On the derived cuboid\jour Canad\. 
Math\. Bull\.\vol 17\issue 4\pages 575-577\yr 1974
\endref
\ref\myrefno{14}\by Chein~E.~Z.\paper On the derived cuboid of an 
Eulerian triple\jour Canad\. Math\. Bull\.\vol 20\issue 4\yr 1977
\pages 509--510
\endref
\ref\myrefno{15}\by Leech~J.\paper The rational cuboid revisited
\jour Amer\. Math\. Monthly\vol 84\issue 7\pages 518--533\yr 1977
\moreref see also Erratum\jour Amer\. Math\. Monthly\vol 85\page 472
\yr 1978
\endref
\ref\myrefno{16}\by Leech~J.\paper Five tables relating to rational cuboids
\jour Math\. Comp\.\vol 32\yr 1978\pages 657--659
\endref
\ref\myrefno{17}\by Spohn~W.~G.\paper Table of integral cuboids and their 
generators\jour Math\. Comp\.\vol 33\yr 1979\pages 428--429
\endref
\ref\myrefno{18}\by Lagrange~J.\paper Sur le d\'eriv\'e du cuboide 
Eul\'erien\jour Canad\. Math\. Bull\.\vol 22\issue 2\yr 1979\pages 239--241
\endref
\ref\myrefno{19}\by Leech~J.\paper A remark on rational cuboids\jour Canad\. 
Math\. Bull\.\vol 24\issue 3\yr 1981\pages 377--378
\endref
\ref\myrefno{20}\by Korec~I.\paper Nonexistence of small perfect 
rational cuboid\jour Acta Math\. Univ\. Comen\.\vol 42/43\yr 1983
\pages 73--86
\endref
\ref\myrefno{21}\by Korec~I.\paper Nonexistence of small perfect 
rational cuboid II\jour Acta Math\. Univ\. Comen\.\vol 44/45\yr 1984
\pages 39--48
\endref
\ref\myrefno{22}\by Wells~D.~G.\book The Penguin dictionary of curious and 
interesting numbers\publ Penguin publishers\publaddr London\yr 1986
\endref
\ref\myrefno{23}\by Bremner~A., Guy~R.~K.\paper A dozen difficult Diophantine 
dilemmas\jour Amer\. Math\. Monthly\vol 95\issue 1\yr 1988\pages 31--36
\endref
\ref\myrefno{24}\by Bremner~A.\paper The rational cuboid and a quartic surface
\jour Rocky Mountain J\. Math\. \vol 18\issue 1\yr 1988\pages 105--121
\endref
\ref\myrefno{25}\by Colman~W.~J.~A.\paper On certain semiperfect cuboids\jour
Fibonacci Quart.\vol 26\issue 1\yr 1988\pages 54--57\moreref see also\nofrills 
\paper Some observations on the classical cuboid and its parametric solutions
\jour Fibonacci Quart\.\vol 26\issue 4\yr 1988\pages 338--343
\endref
\ref\myrefno{26}\by Korec~I.\paper Lower bounds for perfect rational cuboids 
\jour Math\. Slovaca\vol 42\issue 5\yr 1992\pages 565--582
\endref
\ref\myrefno{27}\by Guy~R.~K.\paper Is there a perfect cuboid? Four squares 
whose sums in pairs are square. Four squares whose differences are square 
\inbook Unsolved Problems in Number Theory, 2nd ed.\pages 173--181\yr 1994
\publ Springer-Verlag\publaddr New York 
\endref
\ref\myrefno{28}\by Rathbun~R.~L., Granlund~T.\paper The integer cuboid table 
with body, edge, and face type of solutions\jour Math\. Comp\.\vol 62\yr 1994
\pages 441--442
\endref
\ref\myrefno{29}\by Van Luijk~R.\book On perfect cuboids, \rm Doctoraalscriptie
\publ Mathematisch Instituut, Universiteit Utrecht\publaddr Utrecht\yr 2000
\endref
\ref\myrefno{30}\by Rathbun~R.~L., Granlund~T.\paper The classical rational 
cuboid table of Maurice Kraitchik\jour Math\. Comp\.\vol 62\yr 1994
\pages 442--443
\endref
\ref\myrefno{31}\by Peterson~B.~E., Jordan~J.~H.\paper Integer hexahedra equivalent 
to perfect boxes\jour Amer\. Math\. Monthly\vol 102\issue 1\yr 1995\pages 41--45
\endref
\ref\myrefno{32}\by Rathbun~R.~L.\paper The rational cuboid table of Maurice 
Kraitchik\jour e-print \myhref{http://arxiv.org/abs/math/0111229}{math.HO/0111229} 
in Electronic Archive \myEarXivlink
\endref
\ref\myrefno{33}\by Hartshorne~R., Van Luijk~R.\paper Non-Euclidean Pythagorean 
triples, a problem of Euler, and rational points on K3 surfaces\publ e-print 
\myhref{http://arxiv.org/abs/math/0606700}{math.NT/0606700} 
in Electronic Archive \myEarXivlink
\endref
\ref\myrefno{34}\by Waldschmidt~M.\paper Open diophantine problems\publ e-print 
\myhref{http://arxiv.org/abs/math/0312440}{math.NT/0312440} 
in Electronic Archive \myEarXivlink
\endref
\ref\myrefno{35}\by Ionascu~E.~J., Luca~F., Stanica~P.\paper Heron triangles 
with two fixed sides\publ e-print \myhref{http://arxiv.org/abs/math/0608185}
{math.NT/0608} \myhref{http://arxiv.org/abs/math/0608185}{185} in Electronic 
Archive \myEarXivlink
\endref
\ref\myrefno{36}\by Ortan~A., Quenneville-Belair~V.\paper Euler's brick
\jour Delta Epsilon, McGill Undergraduate Mathematics Journal\yr 2006\vol 1
\pages 30-33
\endref
\ref\myrefno{37}\by Knill~O.\paper Hunting for Perfect Euler Bricks\jour Harvard
College Math\. Review\yr 2008\vol 2\issue 2\page 102\moreref
see also \myhref{http://www.math.harvard.edu/\volna knill/various/eulercuboid/index.html}
{http:/\negskp/www.math.harvard.edu/\textvolna knill/various/eulercuboid/index.html}
\endref
\ref\myrefno{38}\by Sloan~N.~J.~A\paper Sequences 
\myhref{http://oeis.org/A031173}{A031173}, 
\myhref{http://oeis.org/A031174}{A031174}, and \myhref{http://oeis.org/A031175}
{A031175}\jour On-line encyclopedia of integer sequences\publ OEIS Foundation 
Inc.\publaddr Portland, USA
\endref
\ref\myrefno{39}\by Stoll~M., Testa~D.\paper The surface parametrizing cuboids
\jour e-print \myhref{http://arxiv.org/abs/1009.0388}{arXiv:1009.0388} 
in Electronic Archive \myEarXivlink
\endref
\ref\myrefno{40}\by Sharipov~R.~A.\paper A note on a perfect Euler cuboid.
\jour e-print \myhref{http://arxiv.org/abs/1104.1716}{arXiv:1104.1716} 
in Electronic Archive \myEarXivlink
\endref
\ref\myrefno{41}\by Sharipov~R.~A.\paper Perfect cuboids and irreducible 
polynomials\jour Ufa Mathematical Journal\vol 4, \issue 1\yr 2012\pages 153--160
\moreref see also e-print \myhref{http://arxiv.org/abs/1108.5348}{arXiv:1108.5348} 
in Electronic Archive \myEarXivlink
\endref
\ref\myrefno{42}\by Sharipov~R.~A.\paper A note on the first cuboid conjecture
\jour e-print \myhref{http://arxiv.org/abs/1109.2534}{arXiv:1109.2534} 
in Electronic Archive \myEarXivlink
\endref
\ref\myrefno{43}\by Sharipov~R.~A.\paper A note on the second cuboid conjecture.
Part~\uppercase\expandafter{\romannumeral 1} 
\jour e-print \myhref{http://arxiv.org/abs/1201.1229}{arXiv:1201.1229} 
in Electronic Archive \myEarXivlink
\endref
\ref\myrefno{44}\by Sharipov~R.~A.\paper A note on the third cuboid conjecture.
Part~\uppercase\expandafter{\romannumeral 1} 
\jour e-print \myhref{http://arxiv.org/abs/1203.2567}{arXiv:1203.2567} 
in Electronic Archive \myEarXivlink
\endref
\ref\myrefno{45}\by Sharipov~R.~A.\paper Perfect cuboids and multisymmetric 
polynomials\jour e-print \myhref{http://arxiv.org/abs/1203.2567}
{arXiv:1205.3135} in Electronic Archive \myEarXivlink
\endref
\ref\myrefno{46}\by Sharipov~R.~A.\paper On an ideal of multisymmetric polynomials 
associated with perfect cuboids\jour e-print \myhref{http://arxiv.org/abs/1206.6769}
{arXiv:1206.6769} in Electronic Archive \myEarXivlink
\endref
\ref\myrefno{47}\by Sharipov~R.~A.\paper On the equivalence of cuboid equations and 
their factor equations\jour e-print \myhref{http://arxiv.org/abs/1207.2102}
{arXiv:1207.2102} in Electronic Archive \myEarXivlink
\endref
\ref\myrefno{48}\by Sharipov~R.~A.\paper A biquadratic Diophantine equation associated 
with perfect cuboids\jour e-print \myhref{http://arxiv.org/abs/1207.4081}
{arXiv:1207.4081} in Electronic Archive \myEarXivlink
\endref
\ref\myrefno{49}\by Ramsden~J.~R.\paper A general rational solution of an equation 
associated with perfect cuboids\jour e-print \myhref{http://arxiv.org/abs/1207.5339}
{arXiv:1207.5339} in Electronic Archive \myEarXivlink
\endref
\ref\myrefno{50}\by Ramsden~J.~R., Sharipov~R.~A.\paper Inverse problems associated 
with perfect cuboids\jour e-print \myhref{http://arxiv.org/abs/1207.6764}
{arXiv:1207.6764} in Electronic Archive \myEarXivlink
\endref
\ref\myrefno{51}\by Shl\"afli~L.\paper \"Uber die Resultante eines systems mehrerer 
algebraishen Gleihungen\jour Denkschr\. Kaiserliche Acad\. Wiss\. Math\.-Natur\.
Kl\.\vol 4\yr 1852\moreref reprinted in {\eightcyr\char '074}Gesammelte mathematische
Abhandlungen{\eightcyr\char '076}, Band \uppercase\expandafter{\romannumeral 2}
\pages 9--112\publ Birkh\"auser Verlag\yr 1953
\endref
\ref\myrefno{52}\by Cayley~A.\paper On the symmetric functions of the roots of 
certain systems of two equations\jour Phil\. Trans\. Royal Soc\. London\vol 147
\yr 1857\pages 717--726
\endref
\ref\myrefno{53}\by Junker~F.\paper \"Uber symmetrische Functionen von mehreren 
Ver\"anderlishen\jour Mathematische Annalen\vol 43\pages 225--270 \yr 1893
\endref
\ref\myrefno{54}\by McMahon~P.~A.\paper Memoir on symmetric functions of the
roots of systems of equations\jour Phil\. Trans\. Royal Soc\. London\vol 181
\yr 1890\pages 481--536
\endref
\ref\myrefno{55}\by McMahon~P.~A. \book Combinatory Analysis. 
\rm Vol\.~\uppercase\expandafter{\romannumeral 1} and 
Vol\.~\uppercase\expandafter{\romannumeral 2}\publ Cambridge Univ\. Press
\yr 1915--1916\moreref see also Third ed\.\publ Chelsea Publishing Company
\publaddr New York\yr 1984
\endref
\ref\myrefno{56}\by Noether~E.\paper Der Endlichkeitssats der Invarianten
endlicher Gruppen\jour Mathematische Annalen\vol 77\pages 89--92 \yr 1915
\endref
\ref\myrefno{57}\by Weyl~H.\book The classical groups\publ Princeton Univ\.
Press\publaddr Princeton\yr1939
\endref
\ref\myrefno{58}\by Macdonald~I.~G.\book Symmetric functions and Hall polynomials,
\rm Oxford Mathematical Monographs\publ Clarendon Press\publaddr Oxford\yr 1979 
\endref
\ref\myrefno{59}\by Pedersen~P.\paper Calculating multidimensional symmetric
functions using Jacobi's formula\inbook Proceedings AAECC 9, volume 539 of
Springer Lecture Notes in Computer Science\pages 304--317\yr 1991\publ Springer
\endref
\ref\myrefno{60}\by Milne~P.\paper On the solutions of a set of polynomial equations
\inbook Symbolic and numerical computation for artificial intelligence. Computational 
Mathematics and Applications\eds Donald~B.~R., Kapur~D., Mundy~J.~L.\yr 1992\publ
Academic Press Ltd.\publaddr London\pages 89--101
\endref
\ref\myrefno{61}\by Dalbec~J.\book Geometry and combinatorics of Chow forms
\publ PhD thesis, Cornell University\yr 1995
\endref
\ref\myrefno{62}\by Richman~D.~R.\paper Explicit generators of the invariants of 
finite groups\jour Advances in Math\.\vol 124\issue 1\yr 1996\pages 49--76
\endref
\ref\myrefno{63}\by Stepanov~S.~A.\paper On vector invariants of the symmetric group
\jour Diskretnaya Matematika\vol 8\issue 2\yr 1996\pages 48--62
\endref
\ref\myrefno{64}\by Gonzalez-Vega~L., Trujillo~G.\paper Multivariate Sturm-Habicht 
sequences: real root counting on n-rectangles and triangles\jour Revista Matem\'atica 
Complutense\vol 10\pages 119--130\yr 1997
\endref
\ref\myrefno{65}\by Stepanov~S.~A.\paper On vector invariants of symmetric groups
\jour Diskretnaya Matematika\vol 11\issue 3\yr 1999\pages 4--14
\endref
\ref\myrefno{66}\by Dalbec~J.\paper Multisymmetric functions\jour Beitr\"age zur
Algebra und Geom\.\vol 40\issue 1\yr 1999\pages 27--51
\endref
\ref\myrefno{67}\by Rosas~M.~H.\paper MacMahon symmetric functions, the partition 
lattice, and Young subgroups\jour Journ\. Combin. Theory\vol 96\,A\issue 2\yr 2001
\pages 326--340
\endref
\ref\myrefno{68}\by Vaccarino~F.\paper The ring of  multisymmetric functions
\jour e-print \myhref{http://arxiv.org/abs/math/0205233}{math.RA/0205233} 
in Electronic Archive \myEarXivlink
\endref
\ref\myrefno{69}\by Briand~E.\paper When is the algebra of multisymmetric 
polynomials generated by the elementary multisymmetric polynomials?
\jour Beitr\"age zur Algebra und Geom\.\vol 45 \issue 2\pages 353--368
\yr 2004
\endref
\ref\myrefno{70}\by Rota~G.-C., Stein~J.~A.\paper A problem of Cayley from 1857
and how he could have solved it\jour Linear Algebra and its Applications (special 
issue on determinants and the legacy of Sir Thomas Muir)\vol 411\pages 167--253
\yr 2005
\endref
\ref\myrefno{71}\by Briand~E., Rosas~M.~H.\paper Milne's volume function and vector 
symmetric polynomials\jour Journ. Symbolic Comput. \vol 44\issue 5\yr 2009
\pages 583--590
\endref
\ref\myrefno{72}\by Sharipov~R.~A.\paper On a pair of cubic equations associated 
with perfect cuboids\jour e-print \myhref{http://arxiv.org/abs/1208.0308}
{arXiv:12} \myhref{http://arxiv.org/abs/1208.0308}{08.0308} in Electronic Archive 
\myEarXivlink
\endref
\ref\myrefno{73}\by Connel~I.\book Elliptic curve handbook\publ McGill University
\publaddr Montreal\yr 1999\moreref see \myhref{http://www.math.mcgill.ca/connell/}
{http:/\negskp/www.math} \myhref{http://www.math.mcgill.ca/connell/}
{.mcgill.ca/connell}
\endref
\ref\myrefno{74}\paper\myhref{http://www.claymath.org/millennium/Birch\podcherkivanie
and\podcherkivanie Swinnerton-Dyer\podcherkivanie Conjecture/}
{Birch and Swinnerton-Dyer Conjecture}\jour Millennium Prize Problems\publ Clay 
Mathematics Institute\yr 2000\publaddr Cambridge, Massachusetts, USA
\endref
\ref\myrefno{75}\paper\myhref{http://en.wikipedia.org/wiki/Millennium\podcherkivanie
Prize\podcherkivanie Problems}{Millennium Prize Problems}\jour Wikipedia\publ 
Wikimedia Foundation Inc.\publaddr San Francisco, USA 
\endref
\endRefs
\enddocument
\end